\documentclass[12pt, twoside, reqno]{amsart}
\usepackage[margin=1in]{geometry}
\usepackage[T1]{fontenc}
\usepackage{lmodern}
\usepackage{amsmath,amssymb,amsthm,mathtools}
\usepackage{booktabs,array}
\usepackage{xcolor}
\usepackage[colorlinks=true]{hyperref}
\definecolor{darkblue}{RGB}{0,50,150}
\hypersetup{linkcolor=darkblue, citecolor=darkblue, urlcolor=darkblue}
\usepackage{url}
\usepackage{microtype}

\newtheorem{theorem}{Theorem}[section]
\newtheorem{proposition}[theorem]{Proposition}
\newtheorem{lemma}[theorem]{Lemma}
\newtheorem{corollary}[theorem]{Corollary}
\newtheorem{conjecture}[theorem]{Conjecture}
\theoremstyle{definition}
\newtheorem{definition}[theorem]{Definition}
\newtheorem{example}[theorem]{Example}
\newtheorem{question}[theorem]{Question}
\theoremstyle{remark}

\DeclareMathOperator{\st}{st}
\DeclareMathOperator{\Prb}{\mathbb P}
\newcommand{\Sn}{\mathfrak S}
\newcommand{\E}{\mathbb E}

\newcommand{\rot}{R}
\newcommand{\comp}[1]{\overline{#1}}
\newcommand{\qstar}{\frac{3106970325}{3245660803}}

\title[The rotation strategy wins Penney's game for permutations]{The rotation strategy wins Penney's game for permutations}
\author{Himaghna Roy Choudhury}
\address{Independent Researcher}
\email{himaghna.r1@gmail.com}
\date{}
\subjclass[2020]{Primary 05A05; Secondary 05A15, 60C05}
\keywords{Consecutive permutation patterns, Penney's game, cluster method,
first occurrence, computer-assisted proof}

\begin{document}

\begin{abstract}
In Penney's game for permutations, two players choose distinct patterns
of the same length $k$, and independent random variables with a common
continuous distribution are observed until the relative order of $k$
consecutive observations matches one of the patterns; the player who chose
that pattern wins. Elizalde and Lin conjectured that the second player wins
with probability greater than one half by answering the pattern
$a_1a_2\cdots a_k$ with its rotation $a_ka_1\cdots a_{k-1}$. We prove that
this conjecture holds for every $k\ge3$. We also show that for $k\ge8$ the
rotation wins with probability at least $1/2+0.0053/(k+1)$. The proof
writes the advantage of the second player in terms of clusters, which are
chains of overlapping marked occurrences of the two patterns weighted by
$\pm1$. We show that, for each ordering of the observations, the weights of
its clusters sum to $-1$, $0$, or $1$. This gives an explicit bound on the
contribution of arbitrarily long clusters, which proves the result
analytically for $k\ge8$, and exact rational computations prove it for
$3\le k\le7$.
\end{abstract}

\maketitle

\section{Introduction}\label{sec:intro}

Waiting times for patterns in random sequences are a classical topic in
probability and combinatorics. Penney~\cite{Penney} proposed a game in which
two players each choose a word of length three over $\{H,T\}$ and a fair coin
is tossed until one of the words appears; the player whose word appears first
wins. Conway found a simple rule for the odds of such a game, and
Li~\cite{Li} proved a general version of it with a martingale argument.
Guibas and Odlyzko~\cite{GO} developed the correlation polynomials of words
and showed that the second player can always choose a word that beats the
first player's word. Thus the relation ``beats'' is nontransitive.

Consecutive patterns in permutations form a second line of work. A
permutation $\sigma$ contains a consecutive occurrence of $\pi\in\Sn_k$ if
some $k$ adjacent entries of $\sigma$ are in the same relative order as $\pi$.
Elizalde and Noy~\cite{EN03} began the systematic enumeration of
permutations by consecutive occurrences, and Elizalde~\cite{Survey} surveys
the area. The cluster method of Goulden and Jackson~\cite{GJ} counts words
avoiding a set of factors. Elizalde and Noy~\cite{EN} adapted it to
consecutive patterns, which gives the generating function of
pattern-avoiding permutations in terms of signed clusters.

Elizalde and Lin~\cite{EL} combined the two lines of work in a permutation
version of Penney's game. A sequence $X_1,X_2,\ldots$ of independent random
variables with a common continuous distribution is observed, and the first
player whose pattern occurs consecutively wins. They computed all winning
probabilities for patterns of length three, studied pairs of length four,
and found that the game is again nontransitive. For patterns of length four
they estimated the probabilities numerically from finite horizons. They
also proposed a response for the second player at every length.

In this paper we study that response. For a sequence $x_1,\ldots,x_m$ of
distinct real numbers, its \emph{standardization} $\st(x_1\cdots x_m)$ is the
permutation obtained by replacing each entry by its rank. For $\pi\in\Sn_k$,
define
\[
 T_\pi=\min\{n\ge k:
 \st(X_{n-k+1}\cdots X_n)=\pi\}.
\]
For $A=a_1\cdots a_k\in\Sn_k$, put
\[
 B=\rot(A)=a_ka_1\cdots a_{k-1},
\]
the permutation obtained by moving the last entry of $A$ to the front. The
question is whether $\Prb(T_B<T_A)>1/2$ for every $A$.

\begin{example}\label{ex:123}
Let $k=3$, $A=123$ and $B=312$. Each ordering of the first $n$ observations
has probability $1/n!$. At $n=3$ the game ends with $A$ for the ordering
$123$ and with $B$ for the ordering $312$. At $n=4$ it ends with $A$ only
for the ordering $2134$. It ends with $B$ for the four orderings
\[
 1423,\qquad 2413,\qquad 3412,\qquad 4312.
\]
The orderings $1234$, $3124$ and $4123$ also end with $A$, but the first
three observations of each already form $123$ or $312$. Therefore
\[
 \Prb(T_B<T_A,\ T_B\le4)-\Prb(T_A<T_B,\ T_A\le4)
 =\Bigl(\frac16+\frac4{24}\Bigr)-\Bigl(\frac16+\frac1{24}\Bigr)
 =\frac18.
\]
Elizalde and Lin~\cite{EL} computed the full winning probability
$\Prb(T_A<T_B)\approx0.342$, so $B$ wins with probability about $0.658$.
\end{example}

Elizalde and Lin write $\Prb(\sigma\prec\tau)$ for the probability that
$\sigma$ occurs before $\tau$. Their Conjecture~6.1 reads as follows.

\begin{conjecture}[Elizalde and Lin~\cite{EL}]\label{conj:EL}
For every $k\ge3$ and every $\sigma\in\Sn_k$, the permutation
$\tau=\sigma_k\sigma_1\cdots\sigma_{k-1}$ satisfies
$\Prb(\sigma\prec\tau)<1/2$.
\end{conjecture}

For words, the analysis of Guibas and Odlyzko gives the second player's
best response from the correlations of the two words alone. No comparable
formula is known for permutations, since an overlap between two
occurrences constrains only the relative order of the shared entries.
Example~\ref{ex:123} shows the typical difficulty. The advantage of $B$
among early stopping times is easy to compute, but the full race depends
on arbitrarily long chains of overlapping occurrences. Elizalde and Lin
verified Conjecture~\ref{conj:EL} for $k=3$ and supported it numerically
for $k=4$.

In this paper, we prove Conjecture~\ref{conj:EL}.

\begin{theorem}\label{thm:main}
For every $k\ge3$ and every $A\in\Sn_k$,
\[
 \Prb(T_{\rot(A)}<T_A)>\frac12.
\]
\end{theorem}

The result gives a response with a strict advantage for every choice of the
first player. It does not assert that the response is optimal. Elizalde and
Lin give examples in which a different response has a larger winning
probability; for instance, $2134$ beats $2341$ by more than $1234$ does.
Their Problem~6.2 asks for an optimal response.

The proof is organized as follows. In Section~\ref{sec:clusters}, we write
the advantage of $B$ as $F(1)D$, where $F(1)>0$ and $D$ is a difference of
two signed cluster series. In Section~\ref{sec:tails}, we prove that the
signed clusters on a single ordering cancel down to $-1$, $0$ or $1$, and we
use this to bound the infinite tail of $D$. In Section~\ref{sec:overlaps}, we
count the orderings with two overlapping occurrences and show that the first
nonzero coefficient of $D$ favors $B$. Section~\ref{sec:uniform} proves
$D>0$ for every $k\ge8$. Section~\ref{sec:finite} treats $3\le k\le7$ by
exact computation, and Section~\ref{sec:questions} states two open
questions. The computed data are listed in Appendix~\ref{app:tables}.

All occurrences in this paper are consecutive. Each ordering of the first
$n$ observations has probability $1/n!$, so every first-occurrence
probability is independent of the common continuous law. We assume
throughout that $k\ge3$.

\section{Signed clusters and first occurrences}\label{sec:clusters}

We now express the difference of the two winning probabilities through
signed clusters. Fix distinct patterns $A,B\in\Sn_k$. The rotation relation
is not required in this section. Number window starts from zero, and let
\[
 E_i=\{\st(X_{i+1}\cdots X_{i+k})\in\{A,B\}\},
 \qquad p=\Prb(E_i)=\frac2{k!}.
\]
A start $i$ with $E_i$ is a \emph{target occurrence}, or simply an
occurrence.

\begin{definition}\label{def:cluster}
Let $\sigma\in\Sn_n$. A \emph{marked occurrence} of $\sigma$ is a target
occurrence whose start has been selected. A \emph{cluster} of length $n$ is
a permutation $\sigma\in\Sn_n$ together with a set of marked occurrences
such that the first marked start is zero, the last marked start is $n-k$,
and successive marked starts differ by less than $k$. A cluster with $s$
marks has \emph{weight} $(-1)^s$.
\end{definition}

For instance, let $A=123$ and $B=312$. The permutation $4123$ has the
occurrence $312$ at start zero and the occurrence $123$ at start one.
Marking both starts gives a cluster of length four and weight $(-1)^2=1$
whose last marked occurrence has pattern $A$. Marking only start one does not
give a cluster, because the first mark of a cluster is at start zero.

The patterns are distinct and have the same length, so a marked start has
exactly one pattern label. Let $c_n^A$ and $c_n^B$ be the sums of the weights
of clusters of length $n$ whose last marked occurrence has pattern $A$ and
$B$, respectively. Define
\begin{equation}\label{eq:C}
 C_A(z)=\sum_{n\ge k}c_n^A\frac{z^n}{n!},\qquad
 C_B(z)=\sum_{n\ge k}c_n^B\frac{z^n}{n!}.
\end{equation}

Let $\alpha_n$ count permutations of length $n$ avoiding both targets, with
$\alpha_0=1$. Let $u_n^A$ and $u_n^B$ count permutations of length $n$ whose
only target occurrence is the final window, with the indicated pattern. Put
\[
 F(z)=\sum_{n\ge0}\alpha_n\frac{z^n}{n!},\qquad
 U_A(z)=\sum_{n\ge k}u_n^A\frac{z^n}{n!},\qquad
 U_B(z)=\sum_{n\ge k}u_n^B\frac{z^n}{n!}.
\]

\begin{lemma}\label{lem:factorization}
For every $n\ge k$,
\begin{equation}\label{eq:convolution}
 u_n^A=-\sum_{j=k}^{n}\binom nj\alpha_{n-j}\,c_j^A,\qquad
 u_n^B=-\sum_{j=k}^{n}\binom nj\alpha_{n-j}\,c_j^B.
\end{equation}
Equivalently, as formal power series,
\begin{equation}\label{eq:factorization}
 U_A(z)=-F(z)C_A(z),\qquad U_B(z)=-F(z)C_B(z).
\end{equation}
\end{lemma}

\begin{proof}
We prove the first identity; the second follows by exchanging the labels.
Consider pairs $(\sigma,M)$, where $\sigma\in\Sn_n$ ends with $A$ and $M$ is
a set of target occurrences of $\sigma$ that contains the final one. Give
$(\sigma,M)$ the weight $(-1)^{|M|-1}$. For a fixed $\sigma$ with $r$
earlier target occurrences, the weights sum to $(1-1)^r$, which is one if
$r=0$ and zero otherwise. Hence the total weight of all pairs is $u_n^A$.

Now decompose a pair $(\sigma,M)$. List the marked starts in increasing
order, and let $s$ be the smallest marked start such that all successive
marked starts from $s$ onward differ by less than $k$. Every marked start
before $s$ is at most $s-k$, so its window lies among the first $s$
positions. Let $\pi=\st(\sigma_1\cdots\sigma_s)$ carry the marks before
$s$, and let $\gamma=\st(\sigma_{s+1}\cdots\sigma_n)$ carry the remaining
marks, shifted by $s$. Each mark of $\pi$ or $\gamma$ is a target occurrence,
because standardization preserves the relative order inside a window. The
marked permutation $\gamma$ is a cluster of length $j=n-s$ with last pattern
$A$, and $\pi$ is an arbitrary permutation of length $n-j$ with an arbitrary
set of marked occurrences.

Conversely, choose a set of $n-j$ values in $\{1,\ldots,n\}$, a marked
permutation $\pi$ of length $n-j$, and a cluster $\gamma$ of length $j$.
Arranging the chosen values in the order of $\pi$, followed by the remaining
values in the order of $\gamma$, recovers a unique pair $(\sigma,M)$. The
marks of $\pi$ end at least $k$ starts before the first mark of $\gamma$,
so the decomposition of this pair returns $\pi$ and $\gamma$. The weight
factors as
\[
 (-1)^{|M|-1}=-(-1)^{|M_\pi|}\,(-1)^{|M_\gamma|}.
\]
For a fixed permutation $\pi$, the weights $(-1)^{|M_\pi|}$ sum to one if
$\pi$ avoids both targets and to zero otherwise, so their total over all
$\pi$ of length $n-j$ is $\alpha_{n-j}$. Summing over $j$ gives
\eqref{eq:convolution}. Multiplying by $z^n/n!$ and summing over $n$
gives~\eqref{eq:factorization}.
\end{proof}

\section{Cancellation and bounds for the infinite tail}\label{sec:tails}

In this section we prove that the clusters on one ordering cancel almost
completely, and we use this to bound the tail of the cluster series. Fix an
ordering $\sigma\in\Sn_{k+d}$ whose first and last windows are target
occurrences, and let its target starts be
\[
 0=t_1<t_2<\cdots<t_r=d.
\]
For $1\le j\le r$, let $w_j$ be the total weight of the marked sets that
contain $t_1$ and $t_j$, are contained in $\{t_1,\ldots,t_j\}$, and have
successive differences less than $k$. Thus $w_r$ is the total weight of all
clusters on $\sigma$.

\begin{lemma}\label{lem:cancellation}
For each ordering, $w_r\in\{-1,0,1\}$. Moreover, $w_r=0$ unless the actual
target starts connect $0$ to $d$ with successive differences less than $k$.
\end{lemma}

\begin{proof}
Classifying a marked set by its penultimate mark $t_i$ gives
\[
 w_1=-1,\qquad
 w_j=-\sum_{\substack{i<j\\t_j-t_i<k}}w_i\quad(j>1).
\]
Set $S_0=0$ and $S_j=w_1+\cdots+w_j$. Let $\ell(j)$ be the largest index with
$t_{\ell(j)}\le t_j-k$, or zero if there is none. The indices in the sum are
exactly $\ell(j)<i<j$, so $w_j=-(S_{j-1}-S_{\ell(j)})$. Therefore
\[
 S_1=-1,\qquad S_j=S_{j-1}+w_j=S_{\ell(j)}\quad(j>1).
\]
Since $\ell(j)<j$, induction on $j$ gives $S_j\in\{0,-1\}$ for every $j$.
Thus $w_j=S_j-S_{j-1}\in\{-1,0,1\}$. If two consecutive actual starts differ
by at least $k$, then no admissible marked set crosses that gap, so $w_r=0$.
\end{proof}

Write $d_n=c_n^A-c_n^B$. For $d\ge1$, let $G_d$ be the event that targets
occur at starts zero and $d$ and that the actual occurrences between them
have successive differences less than $k$. Every ordering of length $k+d$
contributes $\pm w_r$ to $d_{k+d}$, with the sign given by its last pattern.
By Lemma~\ref{lem:cancellation}, this contribution is zero outside $G_d$ and
has absolute value at most one on $G_d$. Hence
\begin{equation}\label{eq:coefficient-bound}
 \frac{|d_{k+d}|}{(k+d)!}\le\Prb(G_d).
\end{equation}
The same bound applies to $|c_{k+d}^A|/(k+d)!$ and
$|c_{k+d}^B|/(k+d)!$ separately.

Put
\[
 L=2k-2,\qquad q=\frac{2(k-1)}{k!}.
\]
Since $2(k-1)<k!$ for $k\ge3$, we have $q<1$. A block of $k-1$ consecutive
starts uses $L$ observations, and by the union bound the probability that it
contains a target is at most $(k-1)p=q$. The next fact makes connectivity
usable.

\begin{lemma}\label{lem:block}
Suppose that the actual target starts connect $0$ to $d$ with successive
differences less than $k$. Then every block $\{a,a+1,\ldots,a+k-2\}$ with
$0<a$ and $a+k-2<d$ contains a target start.
\end{lemma}

\begin{proof}
Let $x$ be the smallest target start of the chain with $x\ge a$; it exists
because $d\ge a$. Its predecessor $y$ in the chain satisfies $y<a$ and
$x-y<k$. Hence $x\le a+k-2$.
\end{proof}

\begin{lemma}\label{lem:tail}
The two cluster series converge absolutely at one, and
\begin{equation}\label{eq:tail}
 \sum_{d\ge k}\frac{|d_{k+d}|}{(k+d)!}
 \le\frac{8(k-1)}{k!\,[k!-2(k-1)]}.
\end{equation}
\end{lemma}

\begin{proof}
Fix $d\ge k$. The windows at starts zero and $d$ are disjoint. Between them
place blocks of $L$ observations beginning at observations
$k,k+L,k+2L,\ldots$, and keep only the blocks that end before observation
$d$. Their number is $m=\lfloor(d-k)/L\rfloor$. In each block consider the
first $k-1$ window starts; all of these windows lie inside the block. By
Lemma~\ref{lem:block}, the event $G_d$ requires a target among these starts
in every block. The two endpoint windows and the $m$ blocks involve disjoint
sets of observations, so independence gives $\Prb(G_d)\le p^2q^m$.

Each value of $m$ arises for exactly $L$ values of $d\ge k$. Therefore
\[
 \sum_{d\ge k}\Prb(G_d)
 \le p^2L\sum_{m\ge0}q^m
 =\frac{p^2L}{1-q}
 =\frac{8(k-1)}{k!\,[k!-2(k-1)]}.
\]
Equation~\eqref{eq:coefficient-bound} proves~\eqref{eq:tail}, and the same
argument applied to $c_{k+d}^A$ and $c_{k+d}^B$ proves absolute convergence.
\end{proof}

\begin{proposition}\label{prop:advantage}
Let $D=C_A(1)-C_B(1)$. Then
\begin{equation}\label{eq:advantage}
 \Prb(T_B<T_A)-\Prb(T_A<T_B)=F(1)D,
 \qquad 1\le F(1)<\infty.
\end{equation}
\end{proposition}

\begin{proof}
The windows at starts $0,k,2k,\ldots$ are disjoint, so
\[
 \frac{\alpha_n}{n!}\le(1-p)^{\lfloor n/k\rfloor}.
\]
Hence $F(1)$ is finite, and $F(1)\ge\alpha_0=1$. The same estimate shows that
a target appears almost surely. The event that the game ends at time $n$
with pattern $A$ has probability $u_n^A/n!$, so $U_A(1)=\Prb(T_A<T_B)$ and
$U_B(1)=\Prb(T_B<T_A)$. The series $F(1)$, $C_A(1)$ and $C_B(1)$ converge
absolutely by Lemma~\ref{lem:tail}, so by Mertens' theorem the Cauchy
products in~\eqref{eq:convolution} converge to the products of the sums.
Evaluating~\eqref{eq:factorization} at one gives~\eqref{eq:advantage}.
\end{proof}

The following refinement gives finite certificates at the smallest lengths.
For $h\ge1$, let $Q_h$ be the event that a target occurs at start zero and
that its connected chain of actual occurrences reaches a start at least $h$.
The first such start is at most $h+k-2$ by the argument of
Lemma~\ref{lem:block}. Thus $Q_h$ depends only on the first $R=h+2k-2$
observations. Write $r_h=\Prb(Q_h)$.

\begin{lemma}\label{lem:refined-tail}
For $h\ge1$ and $R=h+2k-2$,
\begin{equation}\label{eq:refined-tail}
 \sum_{d\ge R}\frac{|d_{k+d}|}{(k+d)!}
 \le\frac{4(k-1)r_h}{k!-2(k-1)}.
\end{equation}
Consequently,
\begin{equation}\label{eq:certificate}
 D\ge\sum_{d=1}^{R-1}\frac{d_{k+d}}{(k+d)!}
       -\frac{4(k-1)r_h}{k!-2(k-1)}.
\end{equation}
\end{lemma}

\begin{proof}
For $d\ge R$, the event $G_d$ implies $Q_h$, a target at start $d$, and a
target among the first $k-1$ starts of every block of $L$ observations
beginning at observations $R,R+L,\ldots$ and ending before observation $d$.
These events involve disjoint sets of observations. Thus
\[
 \Prb(G_d)\le r_hpq^{\lfloor(d-R)/L\rfloor}.
\]
Summing over $d\ge R$ gives $r_hpL/(1-q)$, which is the right-hand side
of~\eqref{eq:refined-tail}. The single-mark clusters give $c_k^A=c_k^B=-1$,
so $d_k=0$, and~\eqref{eq:certificate} follows.
\end{proof}

\section{Overlap estimates for a rotation}\label{sec:overlaps}

We now count orderings with two overlapping target occurrences, and we show
that the first nonzero coefficient of $D$ is positive while the next ones
are small. For $1\le d<k$ and $U,V\in\Sn_k$, let $N_d(U,V)$ be the number of
permutations of length $k+d$ with pattern $U$ at start zero and pattern $V$
at start $d$.

The two windows share the $k-d$ positions $d+1,\ldots,k$. The orders of $U$
and $V$ on these positions are \emph{compatible} if
$\st(u_{d+1}\cdots u_k)=\st(v_1\cdots v_{k-d})$. In that case the shared
entries form a chain $c_1<\cdots<c_{k-d}$, which splits the remaining values
into $k-d+1$ \emph{gaps}. Let $\alpha_j$ be the number of entries among the
first $d$ positions that lie in gap $j$, and let $\beta_j$ be the number of
entries among the last $d$ positions that lie in gap $j$. The vector
$(\alpha_j)$ is determined by $U$: $\alpha_j$ is the number of
$i\le d$ such that $u_i$ lies between the $j$-th and $(j+1)$-th smallest of
$u_{d+1},\ldots,u_k$. Similarly $(\beta_j)$ is determined by $V$.

\begin{lemma}\label{lem:overlap}
If the orders of $U$ and $V$ on the shared positions are incompatible, then
$N_d(U,V)=0$. Otherwise
\begin{equation}\label{eq:overlap-product}
 N_d(U,V)=\prod_{j=0}^{k-d}\binom{\alpha_j+\beta_j}{\alpha_j},
 \qquad \sum_j\alpha_j=\sum_j\beta_j=d.
\end{equation}
Moreover, $(\alpha_j)=(\beta_j)$ if and only if
$\{u_1,\ldots,u_d\}=\{v_{k-d+1},\ldots,v_k\}$.
\end{lemma}

\begin{proof}
A permutation is determined by the relative order of all its entries. The
first window fixes every comparison among its entries, and the second window
fixes every comparison among its entries. The only remaining comparisons
are between an entry exclusive to the first window and an entry exclusive to
the second. If these two entries lie in different gaps, the shared chain
orders them. Within gap $j$, the $\alpha_j$ entries of the first window
appear in an order fixed by $U$ and the $\beta_j$ entries of the second
window appear in an order fixed by $V$. A relative order of all entries in
the gap is therefore a shuffle of two chains, and there are
$\binom{\alpha_j+\beta_j}{\alpha_j}$ of them. Every choice of shuffles, one
for each gap, gives a total order whose restrictions to the two windows are
$U$ and $V$. This proves~\eqref{eq:overlap-product}. If the shared orders are
incompatible, no total order exists.

For the last statement, let $\gamma=(\gamma_0,\ldots,\gamma_{k-d})$ be a
vector of nonnegative integers with sum $d$. Write a word of length $k$
consisting of $\gamma_0$ letters $x$, one letter $s$, $\gamma_1$ letters $x$,
one letter $s$, and so on, ending with $\gamma_{k-d}$ letters $x$. Let
$\rho(\gamma)\subseteq\{1,\ldots,k\}$ be the set of positions of the letters
$x$. The vector $\gamma$ is recovered from $\rho(\gamma)$ by counting the
letters $x$ between consecutive letters $s$, so $\rho$ is injective. Listing
the values $1,\ldots,k$ of $U$ in increasing order and marking each as
exclusive or shared produces exactly this word for $\gamma=(\alpha_j)$, so
$\rho(\alpha)=\{u_1,\ldots,u_d\}$. In the same way
$\rho(\beta)=\{v_{k-d+1},\ldots,v_k\}$. Injectivity of $\rho$ proves the
claim.
\end{proof}

The case $d=1$ is especially simple. Each window has one exclusive entry,
which lies in gap $u_1-1$ for the first window and in gap $v_k-1$ for the
second. Thus $N_1(U,V)$ is $2$ when $u_1=v_k$ and $1$ otherwise, provided
the shared orders are compatible.

\begin{lemma}\label{lem:leading}
Let $B=\rot(A)$. Then $d_{k+1}\ge1$.
\end{lemma}

\begin{proof}
For $(U,V)=(B,A)$, both shared orders are $\st(a_1\cdots a_{k-1})$, and
$u_1=a_k=v_k$. Thus $N_1(B,A)=2$. For $(U,V)=(A,B)$, we have $u_1=a_1$ and
$v_k=a_{k-1}$, which differ because $k\ge3$. Hence $N_1(A,B)\le1$.

A pattern $U$ overlaps itself in $k-1$ positions only if it is monotone. In
fact, compatibility of $\st(u_2\cdots u_k)$ with $\st(u_1\cdots u_{k-1})$
says that $u_{i+1}>u_i$ if and only if $u_{i+2}>u_{i+1}$, so all adjacent
comparisons have the same direction. Thus $N_1(B,B)\le1$, with equality only
if $B$ is monotone. Suppose that $B=12\cdots k$. Then $A=23\cdots k1$, and the
last two entries of $a_2\cdots a_k$ decrease while the last two entries of
$b_1\cdots b_{k-1}$ increase, so $N_1(A,B)=0$. The case $B=k\cdots21$ is
symmetric.

An ordering of length $k+1$ with targets at both starts carries exactly one
cluster with two marks, of weight one. Therefore
\[
 d_{k+1}=N_1(A,A)+N_1(B,A)-N_1(A,B)-N_1(B,B).
\]
The first term is nonnegative and the second is $2$. Either $B$ is not
monotone, so $N_1(B,B)=0$ and $N_1(A,B)\le1$, or $B$ is monotone, so
$N_1(B,B)=1$ and $N_1(A,B)=0$. In both cases $d_{k+1}\ge1$.
\end{proof}

For instance, $A=123$ and $B=312$ give $N_1(A,A)=1$, $N_1(B,A)=2$ and
$N_1(A,B)=N_1(B,B)=0$, so $d_4=3$. This agrees with
Example~\ref{ex:123}, where $u_4^B-u_4^A=4-1=3$.

\begin{lemma}\label{lem:shuffle}
Let $B=\rot(A)$. For $2\le d<k$,
\begin{align}
 d_{k+d}&\ge-N_d(A,B)-N_d(B,B),\label{eq:negative}\\
 N_d(A,B)+N_d(B,B)&\le
 \begin{cases}
 \binom{2d}{d},&d<k-1,\\[1mm]
 \frac32\binom{2k-2}{k-1},&d=k-1.
 \end{cases}\label{eq:shuffle-bound}
\end{align}
\end{lemma}

\begin{proof}
When $d<k$, all starts $0,1,\ldots,d$ lie less than $k$ apart, so every set
of marks containing $0$ and $d$ is a cluster. If an ordering has $r$
interior target occurrences, their optional marks contribute $(1-1)^r$.
Thus an ordering contributes zero unless its only occurrences are at the
two endpoints, and then its cluster weight is one. A negative contribution to
$d_{k+d}$ therefore requires last pattern $B$, which
proves~\eqref{eq:negative}.

For~\eqref{eq:shuffle-bound}, put $n_j=\alpha_j+\beta_j$, so that
$\sum_jn_j=2d$. By Vandermonde's identity,
\[
 \binom{2d}{d}=\sum_{\gamma}\prod_j\binom{n_j}{\gamma_j},
\]
where $\gamma$ runs over the vectors with $0\le\gamma_j\le n_j$ and
$\sum_j\gamma_j=d$. The vectors $\gamma=(\alpha_j)$ and $\gamma=(\beta_j)$
both occur, and since $\binom{n_j}{\beta_j}=\binom{n_j}{\alpha_j}$ their
terms both equal $N_d(U,V)$. Hence $N_d(U,V)\le\binom{2d}{d}$ always, and
$N_d(U,V)\le\frac12\binom{2d}{d}$ when $(\alpha_j)\ne(\beta_j)$.

By Lemma~\ref{lem:overlap}, the vectors differ exactly when the rank sets of
the exclusive entries differ. For $(B,B)$ these sets are the values of $B$
on its first $d$ and on its last $d$ positions. The two sets of positions
differ when $d<k$, so the sets of values differ. For $(A,B)$ they are the
values of $A$ on positions $1,\ldots,d$ and on positions $k-d,\ldots,k-1$,
which differ when $d<k-1$. Applying the half bound to both pairs when
$d<k-1$, and only to $(B,B)$ when $d=k-1$, proves~\eqref{eq:shuffle-bound}.
\end{proof}

\section{The uniform estimate}\label{sec:uniform}

In this section we prove that $D>0$ for every $k\ge8$. For $d\ge2$ define
\[
 f_d=\binom{2d}{d}\frac{k!}{(k+d)!},\qquad
 S_k=\sum_{d=2}^{k-2}f_d+\frac32f_{k-1},\qquad
 H_k=\frac{8(k-1)}{k!-2(k-1)}.
\]
Splitting $k!D=\sum_{d\ge1}k!\,d_{k+d}/(k+d)!$ into the term $d=1$, the terms
$2\le d\le k-1$, and the terms $d\ge k$, Lemmas~\ref{lem:leading},
\ref{lem:shuffle} and~\ref{lem:tail} give
\begin{equation}\label{eq:master-bound}
 k!D\ge\frac1{k+1}-S_k-H_k.
\end{equation}

\begin{proposition}\label{prop:large}
For $k\ge8$ and $A\in\Sn_k$,
\begin{equation}\label{eq:uniform-positive}
 D\ge\frac{138690478}{3245660803\,(k+1)k!}>0.
\end{equation}
\end{proposition}

\begin{proof}
For $d\ge1$, the identity $(d+1)\binom{2d+2}{d+1}=2(2d+1)\binom{2d}{d}$ gives
\[
 \frac{f_{d+1}}{f_d}
 =\frac{2(2d+1)}{(d+1)(k+d+1)}.
\]
For $d\ge2$ this is at most $4/(k+3)$, since
$4(d+1)(k+d+1)-2(2d+1)(k+3)=4d(d-1)+2(k-1)\ge0$. Hence
$f_d\le f_2\,(4/(k+3))^{d-2}$ for every $d\ge2$, where
$f_2=6/((k+1)(k+2))$. Writing
$S_k=\sum_{d=2}^{k-1}f_d+\frac12f_{k-1}$ and summing a geometric majorant,
\begin{align*}
 (k+1)S_k
 &\le (k+1)\left(\frac{f_2}{1-4/(k+3)}
       +\frac12f_2\left(\frac4{k+3}\right)^{k-3}\right)\\
 &=\frac{6(k+3)}{(k+2)(k-1)}
       +\frac3{k+2}\left(\frac4{k+3}\right)^{k-3}.
\end{align*}
Both expressions on the last line decrease for $k\ge8$. For the first, the
derivative of $(x+3)/((x+2)(x-1))$ has numerator $-x^2-6x-5<0$. For the
second, the prefactor and the base decrease, while the exponent increases
and the base belongs to $(0,1)$.

The expression
\[
 (k+1)H_k=\frac{8(k^2-1)}{k!-2(k-1)}
\]
also decreases for $k\ge8$. Indeed, its denominator at $k+1$ exceeds $k+1$
times its denominator at $k$, because $2k<2(k+1)(k-1)$, whereas its
numerator increases by the factor $k(k+2)/((k-1)(k+1))<k+1$. Evaluating the
three decreasing bounds at $k=8$ yields
\begin{equation}\label{eq:rational-constant}
 (k+1)(S_k+H_k)
 \le\frac{33}{35}+\frac3{10}\left(\frac4{11}\right)^5
       +\frac{504}{40306}
 =\qstar<1.
\end{equation}
Combining~\eqref{eq:master-bound} and~\eqref{eq:rational-constant} gives
$k!D\ge(1-\qstar)/(k+1)$, which is~\eqref{eq:uniform-positive}.
\end{proof}

\section{Exact certificates for the remaining lengths}\label{sec:finite}

We now prove $D>0$ for $3\le k\le7$. For these lengths the uniform estimate
is not sharp enough, so we compute the relevant coefficients exactly. All
counts below are integers, and all bounds and comparisons are computed in the
rational numbers.

For $A=a_1\cdots a_k$, define $\comp A=(k+1-a_1)\cdots(k+1-a_k)$.
Complementation commutes with rotation and preserves all cluster counts and
first-occurrence probabilities. No permutation of length $k\ge2$ equals its
own complement, so it suffices to examine the $k!/2$ permutations that precede
their complements in lexicographic order.

\subsection{Lengths six and seven}
For $k=6,7$, we use the exact overlap counts in~\eqref{eq:overlap-product}
in place of their upper bounds. Define
\begin{equation}\label{eq:pair-certificate}
 J(A)=\frac{d_{k+1}}{k+1}
 -\sum_{d=2}^{k-1}
  \bigl(N_d(A,B)+N_d(B,B)\bigr)\frac{k!}{(k+d)!}
 -H_k.
\end{equation}
Lemmas~\ref{lem:tail}, \ref{lem:leading} and~\ref{lem:shuffle} give
$k!D\ge J(A)$.

\subsection{Lengths three, four and five}
For $k=3,4,5$, the bound $J(A)$ is negative for $3$, $12$ and $18$ of the
representatives, respectively, and we evaluate~\eqref{eq:certificate}
instead. This requires the coefficients
$d_{k+i}$ for $1\le i\le R-1$ and the probability $r_h$. We compute them by
inserting one observation at a time.

To append an observation to an ordering of length $n$, insert it at a rank
$r\in\{0,\ldots,n\}$ and increase each existing rank at least $r$ by one.
Deleting the final entry and standardizing gives each permutation of length
$n+1$ a unique predecessor of length $n$ and a unique insertion rank.
Starting from the two orderings $A$ and $B$ and repeating this insertion
therefore produces every ordering with a target at start zero exactly once.

The state of an ordering records three things: the ranks of its last $k-1$
entries among all $n$ entries, the distance from its latest target start,
and the last $k$ values of $b_i=-S_i$. Here $S_i$ is the partial sum of
Lemma~\ref{lem:cancellation}, indexed by start position: it is constant
between actual target starts and is zero at negative indices. By the proof
of Lemma~\ref{lem:cancellation}, a target at start $i$ gives
\begin{equation}\label{eq:bit-update}
 b_i=b_{i-k},\qquad w_i=b_{i-1}-b_i,
\end{equation}
and a non-target start gives $b_i=b_{i-1}$ and $w_i=0$. In particular each
$b_i$ is $0$ or $1$. A final pattern $A$ contributes $w_i$ to $d_{k+i}$, and
a final pattern $B$ contributes $-w_i$. The initial orderings $A$ and $B$
have distance zero and $b_0=1$, $b_{-1}=\cdots=b_{1-k}=0$.

The ranks of the last $k-1$ entries determine the pattern of the new window,
and the distance determines whether a new occurrence is connected to start
zero. After $k-1$ consecutive starts without a target, no later start can
reconnect, and the ordering contributes nothing more. Such states are
discarded.

\begin{lemma}\label{lem:algorithm}
The insertion recurrence computes the exact coefficients $d_{k+i}$. If
$M_{i,g}$ counts the retained orderings with a target at start $i$ whose
preceding target is at $i-g$, then
\begin{equation}\label{eq:prefix-probability}
 r_h=\sum_{i=h}^{h+k-2}\ \sum_{g=i-h+1}^{k-1}
             \frac{M_{i,g}}{(k+i)!}.
\end{equation}
\end{lemma}

\begin{proof}
By induction on the length, the multiplicity of each state counts the
orderings with that state exactly. The bit update~\eqref{eq:bit-update} and
the signed contribution of the final pattern follow from
Lemma~\ref{lem:cancellation}. Discarded orderings carry no cluster from
start zero to any later start, so they contribute to no coefficient.

For~\eqref{eq:prefix-probability}, partition $Q_h$ by the first target
start $i\ge h$ of the chain. Its predecessor in the chain is below $h$, so
the gap $g$ satisfies $g>i-h$. Connectivity gives $g<k$, hence
$i\le h+k-2$. These events are disjoint and depend only on the first $k+i$
observations, each ordering of which has probability $1/(k+i)!$.
\end{proof}

\subsection{The certificates}

\begin{proposition}\label{prop:finite}
For $3\le k\le7$ and $A\in\Sn_k$, one has $D>0$.
\end{proposition}

\begin{proof}
For $k=3,4,5$ we evaluated $k!$ times the right-hand side
of~\eqref{eq:certificate} for every complement representative, with
$h=13,4,1$, respectively, and with $r_h$ given
by~\eqref{eq:prefix-probability}. The largest lengths used are $19$, $13$
and $13$. For $k=6,7$ we evaluated $J(A)$ for every complement
representative. Table~\ref{tab:certificates} lists the minimum over each
length, and every one is strictly positive:
\[
 \frac{10651931797}{200559744000},\quad
 \frac{8812049}{129729600},\quad
 \frac{98321}{1235520},\quad
 \frac{5539}{89460},\quad
 \frac{4332253}{51768288}
\]
for $k=3,4,5,6,7$, respectively. Since $k!D$ is at least each certified
bound, $D>0$.
\end{proof}

Tables~\ref{tab:k34} and~\ref{tab:k5} list the certified bounds for every
representative with $k\le5$, together with $d_{k+1}$.
Table~\ref{tab:k67} lists the five smallest values of $J(A)$ for $k=6$ and
$k=7$. For $3\le k\le6$ the smallest bound occurs when $B$ is monotone, where
$d_{k+1}=1$; for $k=7$ the monotone case gives the second smallest bound. The complete computation, over
$3+12+60+360+2520$ representatives, took about one minute on a laptop.

We compared these computations with independent enumerations. For every
pattern of length three and four, the coefficients $d_n$ from the insertion
recurrence agree, through $n=8$, with the signed cluster sums obtained by
listing all permutations of length $n$ and applying the recurrence of
Lemma~\ref{lem:cancellation} to their actual occurrences. For the same
patterns, the identity~\eqref{eq:convolution} agrees with a direct count of
first occurrences through $n=8$. The product
formula~\eqref{eq:overlap-product} agrees with a direct count of
permutations for every pair $U,V\in\Sn_k$ with $k\le5$ and every
$1\le d<k$. No disagreement occurred. The data and the code are available in a
public repository~\cite{CodeRepo}.

\begin{proof}[Proof of Theorem~\ref{thm:main}]
Proposition~\ref{prop:large} gives $D>0$ for $k\ge8$, and
Proposition~\ref{prop:finite} gives $D>0$ for $3\le k\le7$.
Proposition~\ref{prop:advantage} then yields
$\Prb(T_B<T_A)>\Prb(T_A<T_B)$. A target occurs almost surely, and the two
distinct patterns cannot first occur at the same time. The two
probabilities therefore sum to one, so $\Prb(T_B<T_A)>1/2$.
\end{proof}

The proof also bounds the size of the advantage. Let $T=\min(T_A,T_B)$.
Since $\alpha_n/n!=\Prb(T>n)$, we have $F(1)=\E[T]$. For $n\ge k$ the first
$n$ observations contain $n-k+1$ windows, so the union bound gives
$\Prb(T>n)\ge1-(n-k+1)p$. With $M=1/p=k!/2$,
\[
 F(1)\ge k+\sum_{m=1}^{M}\Bigl(1-\frac mM\Bigr)=k+\frac{M-1}2\ge\frac{k!}4.
\]

\begin{corollary}\label{cor:quantitative}
For every $k\ge8$ and every $A\in\Sn_k$,
\[
 \Prb(T_{\rot(A)}<T_A)-\frac12\ge\frac{69345239}{12982643212\,(k+1)}.
\]
\end{corollary}

\begin{proof}
By Proposition~\ref{prop:advantage}, the left-hand side equals $F(1)D/2$.
The bound $F(1)\ge k!/4$ and Proposition~\ref{prop:large} give
$F(1)D/2\ge138690478/(8\cdot3245660803\,(k+1))$, which is the stated bound.
\end{proof}

\section{Further questions}\label{sec:questions}

The proof compares the two players only after summing over all stopping
times. Let $a_n=u_n^A$ and $b_n=u_n^B$ count the orderings of length $n$ in
which the game ends at time $n$ with $A$ and with $B$, respectively. The
argument of Lemma~\ref{lem:leading} gives
$b_{k+1}-a_{k+1}=2+N_1(A,A)-N_1(A,B)-N_1(B,B)\in\{1,2,3\}$, and
$b_k=a_k=1$. We computed $\min_A(b_n-a_n)$ by enumerating all permutations
of length $n$ for every $3\le k<n\le11$. The minimum was positive in every
case (Table~\ref{tab:horizon}). Thus one question that could be of interest
is the following.

\begin{question}\label{q:coefficientwise}
Is $b_n>a_n$ for every $k\ge3$, every $A\in\Sn_k$ and every $n>k$?
\end{question}

An affirmative answer would give a second proof of Theorem~\ref{thm:main}.
The most direct injection does not work. For $A=132$ and $B=213$, the
ordering $1243$ first meets a target at its final window, which is $A$.
Rotating its last three entries gives $1324$, which ends with $B$ but
already begins with $A$.

Corollary~\ref{cor:quantitative} shows that the advantage of the rotation
decays at most like $1/k$. The leading term of $D$ suggests that this order
is attained. When $B$ is monotone, $d_{k+1}=1$, and the first term
of~\eqref{eq:master-bound} is only $1/(k+1)$.

\begin{question}\label{q:size}
Is $\min_{A\in\Sn_k}\Prb(T_{\rot(A)}<T_A)-1/2$ of order $1/k$ as
$k\to\infty$?
\end{question}

Elizalde and Lin observed that $\rot(A)$ is not always the best response;
for $A=2341$ the response $2134$ does better. Their Problem~6.2, to find an
optimal response for every $A$, remains open.

\appendix
\section{Computed data}\label{app:tables}

Table~\ref{tab:certificates} summarizes the certificates of
Proposition~\ref{prop:finite}. Tables~\ref{tab:k34}, \ref{tab:k5}
and~\ref{tab:k67} give the bounds for individual representatives. In these
three tables, decimals are truncated, not rounded, so each displayed value
is still a lower bound for $k!D$. Table~\ref{tab:horizon} gives the
finite-horizon data of Question~\ref{q:coefficientwise}.

\begin{table}[ht]
\centering
\caption{Minimum certified lower bounds for $k!D$. For $k\le5$ the bound is
\eqref{eq:certificate} with parameter $h$; the largest permutation length
used is $k+h+2k-3$. For $k=6,7$ the bound is $J(A)$
from~\eqref{eq:pair-certificate}.}
\label{tab:certificates}
\begin{tabular}{@{}rrrrl@{}}
\toprule
$k$&Representatives&$h$&Largest length&Minimum lower bound for $k!D$\\
\midrule
3&3&13&19&$10651931797/200559744000\approx0.0531$\\
4&12&4&13&$8812049/129729600\approx0.0679$\\
5&60&1&13&$98321/1235520\approx0.0795$\\
6&360&\multicolumn{2}{c}{$J(A)$}&$5539/89460\approx0.0619$\\
7&2520&\multicolumn{2}{c}{$J(A)$}&$4332253/51768288\approx0.0836$\\
\bottomrule
\end{tabular}
\end{table}

\begin{table}[ht]
\centering
\caption{Certified lower bounds for $k!D$ for every complement
representative with $k=3$ ($h=13$) and $k=4$ ($h=4$).}
\label{tab:k34}
\begin{tabular}{@{}llrr@{\qquad}llrr@{}}
\toprule
$A$&$B$&$d_{k+1}$&Bound&$A$&$B$&$d_{k+1}$&Bound\\
\midrule
123&312&3&0.3018&1423&3142&2&0.2219\\
132&213&1&0.0784&1432&2143&2&0.1530\\
213&321&1&0.0531&2134&4213&2&0.2276\\
1234&4123&3&0.3852&2143&3214&2&0.2167\\
1243&3124&2&0.2619&2314&4231&2&0.2251\\
1324&4132&2&0.2612&2341&1234&1&0.0679\\
1342&2134&2&0.3076&2413&3241&2&0.3195\\
&&&&2431&1243&2&0.2817\\
\bottomrule
\end{tabular}
\end{table}

\begin{table}[ht]
\centering
\caption{Certified lower bounds for $k!D$ for every complement
representative with $k=5$ ($h=1$).}
\label{tab:k5}
\small
\begin{tabular}{@{}lrr@{\quad}lrr@{\quad}lrr@{}}
\toprule
$A$&$d_6$&Bound&$A$&$d_6$&Bound&$A$&$d_6$&Bound\\
\midrule
12345&3&0.3921&15342&2&0.2613&24351&2&0.2784\\
12354&2&0.2776&15423&2&0.2482&24513&2&0.2865\\
12435&2&0.2549&15432&2&0.2214&24531&2&0.2677\\
12453&2&0.2863&21345&2&0.2493&25134&2&0.2802\\
12534&2&0.2370&21354&2&0.2517&25143&2&0.2575\\
12543&2&0.2511&21435&2&0.2708&25314&2&0.2174\\
13245&2&0.2623&21453&2&0.2655&25341&2&0.2587\\
13254&2&0.2775&21534&2&0.2471&25413&2&0.2725\\
13425&2&0.2557&21543&2&0.2137&25431&2&0.2677\\
13452&2&0.2957&23145&2&0.2629&31245&2&0.2548\\
13524&2&0.2164&23154&2&0.2522&31254&2&0.2386\\
13542&2&0.2548&23415&2&0.2465&31425&2&0.2836\\
14235&2&0.2361&23451&1&0.0795&31452&2&0.2413\\
14253&1&0.0833&23514&2&0.2458&31524&2&0.2716\\
14325&2&0.2734&23541&2&0.2688&31542&2&0.2907\\
14352&2&0.2487&24135&2&0.2203&32145&2&0.2548\\
14523&2&0.2630&24153&2&0.2372&32154&2&0.2448\\
14532&2&0.2327&24315&2&0.2641&32415&2&0.2743\\
15234&2&0.2661&&&&32451&2&0.2640\\
15243&2&0.2756&&&&32514&2&0.2765\\
15324&2&0.2446&&&&32541&2&0.2478\\
\bottomrule
\end{tabular}
\end{table}

\begin{table}[ht]
\centering
\caption{The five smallest values of $J(A)$ among the complement
representatives for $k=6$ (of $360$) and $k=7$ (of $2520$). The largest
values are about $0.3654$ and $0.3643$, respectively.}
\label{tab:k67}
\begin{tabular}{@{}llll@{}}
\toprule
$A$&$B$&$J(A)$&Decimal\\
\midrule
234561&123456&$5539/89460$&0.0619\\
352416&635241&$67241/357840$&0.1879\\
152634&415263&$33727/178920$&0.1885\\
152643&315264&$36709/178920$&0.2051\\
216543&321654&$18461/89460$&0.2063\\
\midrule
1526374&4152637&$4332253/51768288$&0.0836\\
2345671&1234567&$5140085/51768288$&0.0992\\
1642753&3164275&$21355597/103536576$&0.2062\\
3164275&5316427&$8899597/43140240$&0.2062\\
3625147&7362514&$4553501/21570120$&0.2111\\
\bottomrule
\end{tabular}
\end{table}

\begin{table}[ht]
\centering
\caption{The minimum of $b_n-a_n$ over all $A\in\Sn_k$, by exhaustive
enumeration of the permutations of length $n$. A counterexample to
Question~\ref{q:coefficientwise} would be an entry at most zero; none
occurred.}
\label{tab:horizon}
\begin{tabular}{@{}r|rrrrrrrr@{}}
\toprule
$n\backslash k$&3&4&5&6&7&8&9&10\\
\midrule
4&1&&&&&&&\\
5&2&1&&&&&&\\
6&5&6&1&&&&&\\
7&22&32&7&1&&&&\\
8&102&252&56&8&1&&&\\
9&587&2096&469&72&9&1&&\\
10&4120&19756&4685&720&90&10&1&\\
11&30624&201326&50883&7794&990&110&11&1\\
\bottomrule
\end{tabular}
\end{table}

\end{document}